\documentclass{elsart3-1}

\usepackage[english,francais]{babel}
\usepackage[latin1]{inputenc}  
\usepackage[T1]{fontenc}       
\usepackage{color}
\usepackage{amscd}
\usepackage[psamsfonts]{amssymb}
\usepackage{graphics,epsfig,amssymb,bbm,amsmath}
\newtheorem{theorem}{Theorem}[section]

\newtheorem{e-proposition}[theorem]{Proposition}

\newtheorem{e-definition}[theorem]{Definition\rm}

\newtheorem{theoreme}{Th\'eor\`eme}[section]

\newtheorem{proposition}[theoreme]{Proposition}

\renewcommand{\theequation}{\arabic{equation}}
\def\w{\omega}
\def\W{\Omega}

\def\R{{{\Bbb R}}}
\def\C{{{\Bbb C}}}
\def\BR{{{\Bbb R}}}

\def\aaa{{{\mathcal A}}}

\def\ttt{{\mathcal T}}

\def\ue{{\underline e}}
\def\oe{{\overline e}}
\def\E{{{\Bbb E}}}
\def\dive{{\hbox{div}}}

\def\hhh{{{\mathcal H}}}

\def\indic{{{\mathbbm 1}}}

\def\cim{{\delta}}
\def\EnvG{{\tilde \Delta}}

\def\hhhU{{{\mathcal H}^{G,x_0}}}
\def\aaaU{{{\mathcal A}^{G,x_0}}}
\def\haaaU{{\hat{\mathcal A}^{G,x_0}}}

\def\og{\leavevmode\raise.3ex\hbox{$\scriptscriptstyle\langle\!\langle$~}}
\def\fg{\leavevmode\raise.3ex\hbox{~$\!\scriptscriptstyle\,\rangle\!\rangle$}}

\begin{document}

\begin{frontmatter}


\thanks[label2]{Je tiens \`a remercier N. Enriquez \`a qui cette note doit beaucoup.}
\selectlanguage{english}
\title{Markov chains in a Dirichlet Environment and hypergeometric integrals}

\vspace{-2.6cm}

\selectlanguage{francais}
\title{Cha\^\i nes de Markov
en environnement de Dirichlet et intégrales hypergéométriques}


\selectlanguage{english}
\author[authorlabel1]{Christophe Sabot}
\ead{csabot@umpa.ens-lyon.fr}

\address[authorlabel1]{CNRS, UMPA, ENS Lyon, 46 allée d'Italie, 69007 Lyon}

\begin{abstract}
The aim of this text is to establish some relations between Markov
chains in Dirichlet Environments on directed graphs and certain
hypergeometric integrals associated with a particular arrangement
of hyperplanes. We deduce from these relations and the computation
of the connexion obtained by moving one hyperplane of the
arrangement some new relations on important functionals of the
Markov chain.

\vskip 0.5\baselineskip

\selectlanguage{francais}
\noindent{\bf R\'esum\'e} \vskip 0.5\baselineskip \noindent Le but
de ce texte est d'établir une relation entre les chaînes de
Markov en environnement de Dirichlet sur des graphes orientés, et
certaines intégrales hypergéométriques associées à un arrangement
d'hyperplans. Nous déduisons du calcul de la connexion obtenue en
bougeant un hyperplan des relations nouvelles sur des
fonctionnelles importantes de ces marches.

\end{abstract}
\end{frontmatter}

\selectlanguage{francais}
\section*{Version fran\c{c}aise abr\'eg\'ee}
Les lois de Dirichlet jouent un role spécifique dans le contexte
des marches aléatoires en environnement aléatoire, car la loi
moyennée d'une marche aléatoire en environnement de Dirichlet est
la loi d'une marche renforcée sur les arêtes orientées (cf
\cite{EnriquezS}). De façon plus profonde, le but de cette note
est d'établir un lien entre certaines fonctionnelles importantes
de ces marches, étroitement liées à la transformée de Laplace du
temps d'occupation des arêtes, et des intégrales hypergéométriques
associées à un arrangement d'hyperplans particulier. On calcule
alors la connexion obtenue en faisant bouger un hyperplan; elle
induit un système différentiel satisfait par les fonctionnelles de
notre marche aléatoire. Nous pensons que les propriétés
algébriques de cette connexion pourraient donner beaucoup
d'information sur les marches aléatoires qui nous intéressent.

On se place dans le cadre général suivant: on considère un graphe
orienté $G=(V,E)$ contenant un sommet cimetière $\cim$, dont
aucune arête ne sort. On note $U=V\setminus \{\cim\}$, et on
suppose qu'il existe un chemin simple orienté entre tout point $x$
de $U$ et $\cim$. On suppose aussi qu'on a un point $x_0$ dans
$U$, à partir duquel on peut atteindre tout point de $V$ par un
chemin orienté. On note $\oe$ et $\ue$ la tête et la queue d'une
arête $e$. On se donne des poids positifs $\alpha_e$ sur chacune
des arêtes $e$. On construit alors la chaîne de Markov en
environnement aléatoire sur le graphe $G$ de la façon suivante: en
chaque point $x$ de $U$ on tire des probabilités de sortie
$(p_e)_{\ue =x}$, indépendamment suivant une loi de Dirichlet de
paramêtre $(\alpha_e)_{\ue=x}$. On note $\E^{(\alpha)}$ la moyenne
associée sur les environnements, et $P^{(p_e)}$ la loi de la
chaîne de Markov sur le graphe $G$, partant de $x_0$ et arrêtée en
$\cim$, obtenue à partir des probabilités de transition $(p_e)$.
On note alors $G^p$ la fonction de Green de la marche tuée en
$\cim$, de telle sorte que $z_e=G^p(x_0,\ue)p_e$ est le nombre moyen
de traversées de l'arête $e$ avant d'atteindre $\cim$. Le but de
cette note est d'exprimer certaines fonctionnelles comme
$\E^{(\alpha)}[ e^{-<\lambda, z>} p_T]$, où $(\lambda_e)\in
(\R_+^*)^E$ et $p_T$ est la probabilité de l'arbre couvrant $T$
orienté vers $\cim$ (pour la mesure de probabilité naturellement
associé à la chaîne de Markov $P^{(p_e)}$), en fonction
d'intégrales hypergéométriques associées à un arrangement
d'hyperplan particulier. Cette correspondance est obtenue en
faisant le changement de variables $(p_e)\mapsto (z_e)$.

Décrivons maintenant brièvement l'arrangement considéré. On note
$\dive : \R^E\rightarrow \R^U$ l'opérateur $\dive
(\theta)(x)=\sum_{\ue =x} \theta_e -\sum_{\oe=x} \theta_e$ et
$\hhhU$ l'espace affine $\hhhU=\{(z_e),\;
\dive(z)=\delta_{x_0}\}$. Les hyperplans $H_e=\{z_e=0\}$ forment
un arrangement d'hyperplans $\aaaU=(\hhhU,(H_e)_{e\in E})$. Nous
sommes particulierement intéressés par le domaine
$\Delta=\{z_e>0\}$. Les bases de cet arrangement sont les familles
$(H_{e_1}, \ldots ,H_{\vert E\vert -\vert U\vert})$ où $(e_1,
\ldots ,e_{\vert E\vert -\vert U\vert})$ forment le complémentaire
d'un arbre couvrant $T$ de $E$. On note alors $\w_T={dz_{e_1}\over
z_{e_1}} \wedge \cdots \wedge {dz_{e_{\vert E\vert -\vert
U\vert}}\over z_{e_{\vert E\vert -\vert U\vert}}}$, et nous
considérons les intégrales hypergéométriques $I_{\Delta,
T}^{(\alpha)}(\lambda)$ définies en (1). Les intégrales de ce type
ont été intensivement étudiées (cf \cite{OrlikT} pour un texte de
présentation). En particulier, en bougeant un hyperplan ces
intégrales satisfont un système différentiel (\cite{OrlikT}), que
nous calculons dans ce cas précis, cf théorème 3.1. Nous relions
ces intégrales (en fait les intégrales obtenues à partir d'un
graphe bipartite $\hat G$ simplement construit à partir de $G$) à
la transformée de Laplace de la densité d'occupation des arêtes
$z_e$, pondérée par le poids des arbres couvrants du graphe, cf
théorème 2.1.

\selectlanguage{english}
\section{Introduction}
\label{Intro} Among random environments, Dirichlet environments
play a special role since the annealed law of a random walk in a
Dirichlet environment corresponds to the law of a reinforced
random walk, on oriented edges (cf \cite{EnriquezS}). More deeply,
the aim of this text is to show that random walks in Dirichlet
environment have strong relations with certain hypergeometric
integrals, which have been intensively studied (cf e.g.
\cite{OrlikT} for a review text). As a consequence of these
relations, we are able to write a differential system satisfied by
the Laplace transform of certain important functionals of the
walks (this differential system is related to the Gauss-Manin
connection of the arrangement, cf \cite{OrlikT}, chap 8, 10). We
hope this new direction will be useful to understand the
properties of Random Walks in Dirichlet Environment. The complete
proofs of the results announced in this note will appear soon.

Let us now describe the model of  Markov chains in random
Dirichlet environment on directed graphs. Let $G=(V,E)$ be a
finite directed graph, $V$ is the set of vertices and $E\subset
V\times V$ is the set of edges. We denote by $\ue$ (resp. $\oe$)
the origin (resp. the destination) of the edge $e$, so that
$e=(\ue, \oe)$. For simplicity we assume that there is no loop
(i.e. edge of the type $(x,x)$). We suppose that the set of
vertices can be decomposed in $V=U\sqcup \{\cim\}$, and that a
base point $x_0\in U$ is given, such that: there is no edge with
origin $\cim$; for all $x\in U$ there is a directed simple path
from $x$ to $\cim$;  for any $x$ in $U$, there is a directed
simple path from $x_0$ to $x$, with the following simple
definitions.
\begin{e-definition}
i) A simple path from vertices $x$ to $y$, $x\neq y$, is a set of
edges $\{ e_1, \cdots ,e_k\}$ such that there is a list of
distinct vertices $(x_0=x, x_1, \ldots, x_{k}=y)$ such that for
all $j$, $1\le j\le k$, we have either $e_j=(x_{j-1}, x_j)$ or
$e_j=(x_{j}, x_{j-1})$. The path is directed if for all $j$,
$e_j=(x_{j-1}, x_j)$.

ii) A (simple) cycle is a set of edges $\{ e_1, \ldots ,e_k \}$
such that there is a list of distinct vertices $(x_0, \ldots
,x_{k-1})$ such that for all $j$, $1\le j\le k$, we have
either $e_j=( x_{j-1}, x_j)$ or $e_j=(x_{j}, x_{j-1})$, with the
convention that $x_k=x_0$. The cycle is directed if for all $j$,
$e_j=(x_{j-1}, x_j)$.

iii) A spanning tree is a subset $T$ of edges which contains no
cycle and which contains a simple path between any two vertices
$x$ and $y$. The spanning tree $T$ is directed (towards $\cim$) if
for any vertex $x$ in $U$ it contains a unique edge $e$ with
origin $x$. In this case for any vertex $x$ the (simple) path from
$x$ to $\cim$ in $T$ is directed.
\end{e-definition}

Now, we construct some Markov chain on $G$ killed at $\{\delta\}$.
We define the set of environments
$$
\EnvG= \{(p_e)\in (0,1]^{E}, \;\;\; \forall x\in U, \;\;
\sum_{e,\;\ue=x}p_e=1\}.
$$
An element $(p_e)$ of $\EnvG$ defines a Markov chain on $V$,
starting at $x_0$, stopped at $\cim$, and  with exit probabilities
$(p_e)_{\ue=x}$ at any point $x$ in $U$.

Let us now consider a set of positive weights $ (\alpha_e)_{e\in
E}. $ For all $x$ in $U$ we set $ \beta_x=\sum_{\ue =x} \alpha_e.
$ We endow the set of environments $\EnvG$ with the probability
$\mu^{(\alpha)}$ with density
$$
{\prod_{x\in U} \Gamma(\beta_x)\over  \prod_{e\in E}
\Gamma(\alpha_e)} \left( \prod_{e\in E} p_e^{\alpha_e -1}\right)
d\lambda_{\tilde \Delta},
$$
where $d\lambda_{\tilde \Delta}$ is the measure on $\tilde \Delta$
given by $ d\lambda_{\tilde \Delta} =\prod_{e\in \tilde E} dp_e, $
where $\tilde E$ is obtained from $E$ by removing arbitrarily, for
each vertex $x$, one edge with origin $x$ (one can easily see that
$d\lambda_{\tilde \Delta}$ is independent of this choice).
 This means that the transition probabilities are chosen
independently at each vertex $x$ under a Dirichlet law of
parameter $(\alpha_e)_{\ue=x}$.

We denote by $G^p_U$ the Green function of the Markov chain with
transition probabilities $(p_e)$ killed at $\cim$, i.e.
$$
G_U^p(x,y)= E_x^{(p_e)} \left[ \sum_{k=0}^{T_\cim-1}
\indic_{X_k=y} \right] = (I-P_U)^{-1}_{x,y},
$$
where $T_\cim$ is the first hitting time of $\cim$, and $P_U$ is
the transition matrix of the Markov chain restricted to $U\times
U$. For all edge $e\in E$, $ z_e=G_U^p(x_0,\ue) p_e $ is equal to the
expected number of crossings of the edge $e$ before the killing
time $T_\delta$. These values give of course considerable
information on the behavior of the random walk $P^{(p_e)}_{x_0}$.

\section{The arrangement $(\hhh^{G, x_0}, (H_e)_{e\in E})$, and the change of variables.}

We suppose that we have a graph $G=(V,E)$ and a base point $x_0$
as in section 1, and some (not necessarily positive) weights
$(\alpha_e)_{e\in E}$. We define the divergence operator
$\dive:\R^E\rightarrow \R^U$ by
$$
\dive(\theta)(x)=\sum_{\ue=x} \theta_e -\sum_{\oe=x} \theta_e,
\;\;\; (\theta_e)\in \R^E, \;x\in U,
$$
and $\hhh^{G,x_0}$ the affine space by
$$ \hhh^{G,x_0}=\left\{ (z_e)\in \R^E, \;\;\;
\dive(z)=\delta_{x_0}\right\},
$$
where $\delta_{x_0}$ is the Dirac mass at ${x_0}$. The hyperplanes
$ H_e=\{z_e=0\}\cap \hhhU$ define an arrangement of hyperplanes in
the affine space $\hhhU$ that we denote
$$
\aaa^{G,x_0}=(\hhh^{G,x_0}, (H_e)_{e\in E}).
$$
The complement
$\hhhU\setminus\cup_e H_e$ determines some connected components,
and we will be specially interested in a particular connected
component $ \Delta=\{(z_e)\in \hhhU, \; z_e>0,\; \forall e\in
E\}$. The crucial remark is that the function $(z_e)=(G^p_U(x_0,
\ue)p_e)$ takes its values in $\Delta$ for all environment $(p_e)\in \EnvG$,
thanks to the third assumption on $G,x_0$. In particular $\Delta$
is not empty. We recall that a basis of an arrangement is a maximal
free subfamily of $(H_e)$.
\begin{e-proposition}
The bases of the arrangement $\aaaU$ are exactly the subsets
$\{H_e\}_{e\in T^c}$, for all spanning trees $T$.
\end{e-proposition}
We suppose given on $\hhhU$ an arbitrary orientation. To any
spanning tree $T$, we associate the differential form, with
logarithmic poles
$$
\w_T ={dz_{e_1}\over z_{e_1}} \wedge \cdots \wedge {dz_{e_{\vert
E\vert -\vert U\vert}}\over z_{e_{\vert E\vert -\vert U\vert}}},
$$
where $\{e_1, \ldots ,e_{\vert E\vert -\vert U\vert}\}=T^c$ and is
ordered so that $\w_T$ is positively oriented.

We will be interested in the following integrals (when they are
well-defined): for all $(\lambda_e) \in (\C)^E$,
$\Re(\lambda_e)>0$ we set
$$
 I_{\Delta, T}^{(\alpha)}(\lambda)=\int_\Delta
e^{-<\lambda, z>} \prod_e z_e^{\alpha_e} \w_T. \hspace{40pt} (1)
$$

We want to relate some functional of the Markov chain defined in
section 1 with hypergeometric integrals of the type (1). The
Markov chain in Dirichlet environment defined on $G$ as in section
(1) is not naturally related to the hypergeometric integral (1)
for the graph $G$, but for a closely related graph $\hat G$.

The strategy is to construct a bipartite graph $\hat G=(\hat V,
\hat E)$ from the graph $G=(V,E)$ by duplication of the vertices
of $U$. We define $\hat V=\hat U\cup \{\delta \}$, where $ \hat
U=\{-,+\}\times U.$ To simplify notations, we denote $(-,x)$
(resp. $(+,x)$) by $x_-$ (resp. $x_+$). We construct the set of
edges $\hat E$ as follows: to any edge $e\in E$ we associate an
edge $\hat e\in \hat E$ by $\hat e=(x_+,y_-)$ if $e=(x,y)$ and
$y\neq \delta$, and $\hat e=(x_+, \delta)$ if $y=\delta$. To any
vertex $x$ in $U$ we associate an edge $\hat e_x=(x_-, x_+)$. We
define the set of edges of $\hat G$ by $ \hat E=\{\hat e, \; e\in
E\}\cup\{ \hat e_x, \; x\in U\}. $

Then, we define the affine space $ \hat\hhh^{G,x_0}=\hhh^{\hat G,
(-,x_0)}=\{ (\hat z_e)\in \BR^{\hat E},\;\;\; \dive ((\hat
z_e))=\delta_{(-,x_0)}\}. $ It is clear that the affine spaces
$\hhh^{G,x_0}$ and $\hat \hhh^{G,x_0}$ are isomorphic.
Thus we always write $\hhh^{G,x_0}$ for $\hat \hhh^{G,x_0}$, and
we set $z_{\hat e_x}=\sum_{\ue =x} z_e$ for all vertex $x\in U$.
The arrangement associated with the graph $\hat G$ is then $
\haaaU =(\hhh^{G,x_0}, (H_e)_{e \in \hat E}).$ Compared with the
arrangement associated with $G$, we see that they have the same
underlying affine space $\hhh^{G,x_0}$, but the arrangement of
$\hat G$ is richer since it contains all the hyperplanes of the
arrangement of $G$ plus the hyperplanes $ H_{\hat e_x}=\{z_{\hat
e_x}=0\}=\{\sum_{\ue=x} z_e =0\}, $ for all vertex $x$ in $U$.

Comparing the spanning trees, it is clear that any spanning tree
$T$ of $G$ can be extended into a spanning tree $\hat T$ of $\hat
G$ by adding the edges $\{e_x, \;x\in U\}$. Moreover, if $\hat T$
is a directed spanning tree of $\hat G$, then all the edges $\hat
e_x$ are contained in $\hat T$ (cf definition of section 2), thus
removing the edges $\{\hat e_x, \; x\in U\}$ we get a  directed
spanning tree of $G$. Thus, the directed spanning trees of $G$ and
$\hat G$ are the same.

We give now the following weights to the edges of $\hat E$: $
\alpha_{\hat e}=\alpha_e,\;\; e\in E,\;\; \alpha_{\hat
e_x}=-\beta_x. $ It is clear that the domain $\Delta$ of the
arrangement $\aaaU$ is also a domain of the arrangement $\haaaU$
since $z_{\hat e_x}=\sum_{\ue=x} z_e$ is positive if all the $z_e$
are positive. For all $\lambda \in \C^{\hat E}$,
$\Re(\lambda_e)>0$, and all spanning tree $\hat T$ of $\hat G$, we
denote by
$$
\hat I_{\hat T,\Delta}^{(\alpha)}(\lambda)=\int_\Delta
e^{-<\lambda, z>} \left( \prod_{e\in \hat E} z_e^{\alpha_e}\right)
\w_T.
$$
the integral (1) for the graph $\hat G$ and the weights
$(\alpha_{e}, e\in \hat E)$ such defined. Without loss of
generality we can consider a function $\lambda$ which vanishes on
the edges $\hat e_x$ (indeed, since $\hat z_{\hat e_x}= \sum_{\ue
=x} \hat z_{\hat e}$), so we may consider $\hat I_{\hat T,
\Delta}(\lambda)$ as a function on $\C^E$.
\begin{theorem}
For all directed spanning tree $T$ of $ G$, and for all $\lambda
\in \C^E$, $\Re(\lambda_e)>0$, we have the following equality (and
the integral are well-defined):
$$
C_\alpha  \hat I_{\Delta,\hat
T}^{(\alpha)}(\lambda)=\E^{(\alpha)}\left[ e^{-<\lambda ,z>}
\left( {\prod_{e\in T} p_e \over \det(I-P_U)}\right) \right],
$$
where on the left-hand side 
$C_\alpha=\prod_{x} \Gamma(\beta_x)/ \prod_{e} \Gamma(\alpha_e)$ 
and $\hat T$ is the associated
directed spanning tree in $\hat G$, and on the right-hand side
$z_e= G^p(x_0,\ue)p_e$ is the expected number of visits of the
edge $e$ by the Markov chain on $G$ starting at $x_0$ 
with transition probabilities
$(p_e)$.
\end{theorem}
Sketch of the proof: We make the (not so easy) change of variables
$i:\EnvG\rightarrow \Delta$, given by $(p_e)\mapsto z_e=G^p_U(x_0,
\ue) p_e$. \qed

 This last formula has an important probabilistic
meaning. Indeed, we see that we get the Laplace transform of the
occupation density of the edge, which is an important probabilistic quantity,
weighted by the probability of a directed spanning tree (for the
natural probability measure associated with $(p_e)$). This law on
spanning trees contains a lot of information on the initial Markov
chain. For example, by the Wilson algorithm (cf \cite{Wilson}),
the law of the unique (directed) simple path from $x_0$ to
$\delta$ is the law of the loop erased random walk from $x$ to
$\delta$, obtained from the Markov chain $P^{(p_e)}$. Moreover,
summing on all directed spanning trees we exactly get the Laplace
transform of the occupation density of the edges.

\section{The Gauss-Manin connection}
We consider a graph $G$ and some weights $\alpha_e$ such that the
integrals $I_{\Delta, T}(\lambda)$ are well-defined (actually, we
are mainly interested in the graph $\hat G$ obtained from a graph
$G$ with positive weights, as described in section 2). We denote
by $\ttt$ the set of spanning trees of $G$. To any cycle $C$ of
$G$ (given with an arbitrary orientation) we associate a linear
form $l_C$ on $\R^E$, and to any  simple path $\sigma$ from $x_0$
to $\cim$ we associate a  linear form $l_\sigma$, given by
$$l_C(\lambda)=\sum_{e\in C} \epsilon_C(e) \lambda_e, \;\;\;
l_\sigma (\lambda)=\sum_{e\in \sigma} \epsilon_\sigma(e)
\lambda_e,\;\;\; (\lambda_e)\in \R^E,
$$
where $\epsilon_C(e)$ is equal to $+1$ (resp. -1) if $e$ is
directed according to (resp. in opposition to) the orientation of
$C$, and $\epsilon_\sigma(e)$ is equal to $1$ (or $-1$) if $e$ is
directed from $x_0$ to $\delta$ (or from $\delta$ to $x_0$) in the
path $\sigma$. We define the manifold $ M=\C^E\setminus \cup_{C
\hbox{ cycle}} \ker l_C, $ where the logarithmic differential
forms ${d l_C\over l_C}$ are well-defined.

Let $(f_T)_{T\in \ttt}$ be the canonical base of $\R^\ttt$. To all
simple cycle $C$ of $G$, we associate the endomorphism  $\W_C$ of
$\R^\ttt$, given by
$$
\W_C (f_T)= \sum_{e\in C} \epsilon_C(e_0,e) \alpha_{e}
f_{T\cup\{e_0\}\setminus \{e\}},
$$
if there is an edge $e_0$ in $T^c$ such that $C\subset
T\cup\{e_0\}$ and $\W_C(f_T)=0$ otherwise  (in the last formula,
$\epsilon_C(e_0,e)$ is equal either to $+1$, resp. $-1$, if the
directions of $e$ and $e'$ are the same, resp. opposite, in the
cycle $C$). To any simple path $\sigma$ from $x_0$ to $\cim$, we
associate the diagonal endomorphism $\W_\sigma$ of $\R^\ttt $,
with coefficients $ (\W_\sigma)_{T,T}= 1$ if $\sigma\subset T$,
and $0$ otherwise. It is clear that $\W_\sigma$ is a projector on
the subspace generated by the spanning trees containing $\sigma$.
Concerning the matrices $\W_C$ we easily get $
(\W_C)^2=(\sum_{e\in C} \alpha_e) \W_C, $ which means that
${\W_C\over \sum_{e\in C}\alpha_e}$ is a projector.
\begin{theorem}\label{Gauss}
The vector $ I_{\Delta}(\lambda)=\left( I_{\Delta,
T}(\lambda)\right)_{T\in \ttt}$  satisfies $ d
I_{\Delta}(\lambda)= -\W I_{\Delta}(\lambda), $ on $M\cap
\{\Re(\lambda_e)>0\}$, where $d$ is the differential operator in
the variables $(\lambda)$, and $\W$ is the $\ttt\times \ttt$
matrix (with differential form coefficients), given by
$$
\W=\sum_{\sigma} dl_\sigma (\lambda)  \W_\sigma +\sum_C
{dl_C(\lambda)\over l_C(\lambda)} \W_C,
$$
where the first sum is taken over all simple paths from $x_0$ to
$\delta$ and the second summation is taken over all simple cycles
of $G$.
\end{theorem}
Proof: If $T$ is a spanning tree, and $e$ an edge in $T^c$, there
is a unique cycle $C^{e}_T$ contained in $T\cup\{e\}$, that we
orient by convention according to the direction of $e$. Moreover,
there is a unique simple path $\sigma_T$ from $x_0$ to $\delta$ in
$T$, that we orient from $x_0$ to $\delta$. The proof is based on
two types of relation. The first one is of cohomological nature:
let $e_0$ be in $T^c$ then
$$
l_{C^{e_0}_T}(\lambda)\int_\Delta z_{e_0} e^{-<\lambda ,z>}
\left(\prod_{e\in E} z_e^{\alpha_e}\right) \w_T = \left(
\sum_{e'\in C^{e_0}_T} \epsilon_{C^{e_0}_T}(e') I_{\Delta,
T\cup\{e_0\}\setminus \{e'\}}(\lambda)\right).
$$
The second relation comes from the condition
$\dive(z)=\delta_{x_0}$, which implies that for all
$(\lambda_e)\in \R^E$ and for all spanning tree $T$: $ <z,
\lambda>=l_{\sigma_T}(\lambda)+ \sum_{e\in T^c} z_e
l_{C^e_T}(\lambda)$. \qed

The operator $ \nabla= d+\W$ defines a connection on the vector
bundle $M\times \R^\ttt$, and theorem \ref{Gauss} says that
$I_\Delta(\lambda)$ is a flat section of this connection. This
connection is actually integrable, and we describe the structure
equations. Let us introduce a definition. We define the genus of a
subset $S\subset E$ as the size of the maximal free family of
cycles contained in $S$ (we say that a family $(C_1, \ldots ,C_k)$
of cycles is free if the associated functions $(\chi_{C_1},\ldots
,\chi_{C_{k}})$ are free, where $\chi_C=\sum_C \epsilon^C(e) \delta_e$).
\begin{proposition}
The matrices $\W_C$ and $\W_\sigma$ satisfy the following
commutation relations:

i) $[\W_{C}, \W_{C'}]=0$ if $C$ and $C'$ are disjoint or if the
genus of $C\cup C'$ is not 2.

ii) $[\W_\sigma, \W_{\sigma'}]=0$ for all $\sigma$, $\sigma'$.

iii) $[\W_C, \W_\sigma]=0$ if $C$ and $\sigma$ are disjoint or the
genus of $\sigma\cup C$ is not 1.

iv) If $S\subset C$ has genus 2, then either it contains 2
disjoint cycles, or exactly 3 cycles. In the last case, these 3
cycles $C_1$, $C_2$, $C_3$ satisfy the following commutation
relation
 $[\W_{C_1}+\W_{C_2}+\W_{C_3}, \W_{C_i}]=0$, for all $i=1,\; 2,\;3$.

v) If $\sigma_1$ and $\sigma_2$ are two simple paths from $x_0$ to
$\delta$ such that the genus of $\sigma_1\cup\sigma_2$ is 1, i.e.
such that $\sigma_1\cup\sigma_2$ contains a unique cycle $C$, then
$ [ \W_{\sigma_1}+\W_{\sigma_2}, \W_C]=0$.

These relations imply that the connection $\nabla$ is integrable,
i.e. that $\nabla^2=\W\wedge \W =0$.
\end{proposition}
Sketch of the proof: the commutation relations i), ii), iii), iv),
v) imply that $\W\wedge \W=0$, once we remark that $
{dl_{C_1}\over l_{C_1}}\wedge {dl_{C_2}\over l_{C_2}}-
{dl_{C_1}\over l_{C_1}}\wedge {dl_{C_3}\over l_{C_3}}
+{dl_{C_2}\over l_{C_2}}\wedge {dl_{C_3}\over l_{C_3}}=0, $ for
$C_1,\; C_2,\; C_3$ in the configuration iv), and that
$dl_{\sigma_1}\wedge {dl_C\over l_C}=dl_{\sigma_2}\wedge
{dl_C\over l_C}$ in the configuration v). The commutation
relations are obtained by direct computation. \qed




\end{document}